\documentclass{amsart}
\usepackage{CJKutf8}
\usepackage{units}
\usepackage{amstext}
\usepackage{amsthm}
\usepackage{amssymb}
\usepackage{esint}

\makeatletter
\numberwithin{equation}{section}
\numberwithin{figure}{section}
\theoremstyle{plain}
\newtheorem{thm}{Theorem}

\makeatother

\allowdisplaybreaks
\usepackage[sort&compress,square,numbers]{natbib}
\global\long\def\relphantom#1{\mathrel{\phantom{{#1}}}}
\global\long\def\Zp{\mathbb{Z}_{p}}

\begin{document}
\title[Identities for the generalized degenerate Euler polynomials] {Identities of symmetry for the generalized degenerate Euler polynomials}

\author{Dae San Kim}
\address{Department of Mathematics, Sogang University, Seoul 121-742, Republic
	of Korea}
\email{dskim@sogang.ac.kr}

\author{Taekyun Kim}
\address{Department of Mathematics, Kwangwoon University, Seoul 139-701, Republic
	of Korea}
\email{tkkim@kw.ac.kr}

\keywords{Identities of symmetry, Generalized degenerate Euler polynomial, fermionic $p$-adic integral}

\subjclass[2010]{11B68, 11B83, 11C08, 65D20, 65Q30, 65R20}

\begin{abstract}
In this paper, we give some identities of symmetry for the generalized
degenerate Euler polynomials attached to $\chi$ which are derived
from the symmetric properties for certain fermionic $p$-adic integrals
on $\mathbb{Z}_{p}$.
\end{abstract}

\maketitle
\section{Introduction and preliminaries}

Let $p$ be a fixed odd prime. Throughout this paper, $\mathbb{Z}_{p},\mathbb{Q}_{p}$
and $\mathbb{C}_{p}$ will be the ring of $p$-adic integers, the
field of $p$-adic rational numbers and the completion of the algebraic
closure of $\mathbb{Q}_{p}$, respectively. 

The $p$-adic norm $\left|\cdot\right|_{p}$ in $\mathbb{C}_{p}$
is normalized as $\left|p\right|_{p}=\frac{1}{p}$. Let $f\left(x\right)$
be continuous function on $\mathbb{Z}_{p}$. Then the fermionic $p$-adic
integral on $\mathbb{Z}_{p}$ is defined as 
\begin{align}
I_{-1}\left(f\right) & =\int_{\Zp}f\left(x\right)d\mu_{-1}\left(x\right)\label{eq:1}\\
 & =\lim_{N\rightarrow\infty}\sum_{x=0}^{p^{N}-1}f\left(x\right)\left(-1\right)^{x},\quad\left(\text{see \cite{key-10}}\right).\nonumber 
\end{align}

From (\ref{eq:1}), we note that 
\begin{equation}
I_{-1}\left(f_{n}\right)+\left(-1\right)^{n-1}I_{-1}\left(f\right)=2\sum_{l=0}^{n-1}\left(-1\right)^{n-1-l}f\left(l\right),\quad\left(\text{see \cite{key-7}}\right),\label{eq:2}
\end{equation}
where $n\in\mathbb{N}$. 

As is well known, the Euler polynomials are defined by the generating function
\begin{equation}
\int_{\Zp}e^{\left(x+y\right)t}d\mu_{-1}\left(y\right)=\frac{2}{e^{t}+1}e^{xt}=\sum_{n=0}^{\infty}E_{n}\left(x\right)\frac{t^{n}}{n!}.\label{eq:3}
\end{equation}

When $x=0$, $E_{n}=E_{n}\left(0\right)$ are called the Euler numbers
(see \cite{key-1,key-2,key-3,key-4,key-5,key-6,key-7,key-8,key-9,key-10,key-11,key-12,key-13,key-14,key-15,key-16,key-17,key-18,key-19}).

For a fixed odd integer $d$ with $\left(p,d\right)=1$, we set 
\begin{align*}
X & =\lim_{\underset{N}{\leftarrow}}\nicefrac{\mathbb{Z}}{dp^{N}\mathbb{Z}},\quad X^{*}=\bigcup_{\substack{0<a<dp\\
\left(a,p\right)=1
}
}\left(a+dp\mathbb{Z}_{p}\right),\\
a+dp^{N}\mathbb{Z}_{p} & =\left\{ x\in X\mid x\equiv a\pmod{dp^{N}}\right\}, 
\end{align*}
where $a\in\mathbb{Z}$ lies in $0\le a<dp^{N}$. 

It is known that 
\[
\int_{\Zp}f\left(x\right)d\mu_{-1}\left(x\right)=\int_{X}f\left(x\right)d\mu_{-1}\left(x\right),\quad\left(\text{see \cite{key-7,key-8,key-10}}\right),
\]
where $f$ is a continuous function on $\mathbb{Z}_{p}$. 

Let $d\in\mathbb{N}$ with $d\equiv1\pmod{2}$ and let $\chi$ be
a Dirichlet character with conductor $d$. Then the generalized
Euler polynomials attached to $\chi$ are defined by the generating
function 
\begin{equation}
\left(\frac{2\sum_{a=0}^{d-1}\left(-1\right)^{a}\chi\left(a\right)e^{at}}{e^{dt}+1}\right)e^{xt}=\sum_{n=0}^{\infty}E_{n,\chi}\left(x\right)\frac{t^{n}}{n!}.\label{eq:4}
\end{equation}
In particular, for $x=0$, $E_{n,\chi}=E_{n,\chi}\left(0\right)$ are
called the generalized Euler numbers attached to $\chi$. 

For $d\in\mathbb{N}$ with $d\equiv1\pmod{2}$, by (\ref{eq:2}),
we get 
\begin{align}
 & \int_{X}\chi\left(y\right)e^{\left(x+y\right)t}d\mu_{-1}\left(y\right)\label{eq:5}\\
 & =\frac{2\sum_{a=0}^{d-1}\left(-1\right)^{a}\chi\left(a\right)e^{at}}{e^{dt}+1}e^{xt}\nonumber \\
 & =\sum_{n=0}^{\infty}E_{n,\chi}\left(x\right)\frac{t^{n}}{n!},\quad\left(\text{see \cite{key-9,key-10,key-11}}\right).\nonumber 
\end{align}

From (\ref{eq:5}), we have 
\begin{equation}
\int_{X}\chi\left(y\right)\left(x+y\right)^{n}d\mu_{-1}\left(y\right)=E_{n,\chi}\left(x\right),\quad\left(n\ge0\right).\label{eq:6}
\end{equation}

Carlitz considered the degenerate Euler polynomials given by
the generating function 
\begin{align}
 & \frac{2}{\left(1+\lambda t\right)^{\frac{1}{\lambda}}+1}\left(1+\lambda t\right)^{\frac{x}{\lambda}}\label{eq:7}\\
 & =\sum_{n=0}^{\infty}\mathcal{E}_{n}\left(x\mid\lambda\right)\frac{t^{n}}{n!},\quad\left(\text{see \cite{key-3}}\right).\nonumber 
\end{align}

Note that $\lim_{\lambda\rightarrow0}\mathcal{E}_{n}\left(x\mid\lambda\right)=E_{n}\left(x\right)$,
$\left(n\ge0\right)$. 

From (\ref{eq:2}), we note that 
\begin{align}
 & \int_{X}\left(1+\lambda t\right)^{\frac{x+y}{\lambda}}d\mu_{-1}\left(y\right)\label{eq:8}\\
 & =\frac{2}{\left(1+\lambda t\right)^{\frac{1}{\lambda}}+1}\left(1+\lambda t\right)^{\frac{x}{\lambda}}\nonumber \\
 & =\sum_{n=0}^{\infty}\mathcal{E}_{n}\left(x\mid\lambda\right)\frac{t^{n}}{n!}.\nonumber 
\end{align}

Thus, by (\ref{eq:8}), we get 
\begin{equation}
\int_{X}\left(y+x\mid\lambda\right)_{n}d\mu_{-1}\left(y\right)=\mathcal{E}_{n}\left(x\mid\lambda\right),\quad\left(n\ge0\right),\label{eq:9}
\end{equation}
where $\left(x\mid\lambda\right)_{n}=x\left(x-\lambda\right)\cdots\left(x-\left(n-1\right)\lambda\right)$,
for $n\geq1$, and $\left(x\mid\lambda\right)_{0}=1$. 

From (\ref{eq:2}), we can derive the following equation:

\begin{align}
 & \int_{X}\chi\left(y\right)\left(1+\lambda t\right)^{\frac{x+y}{\lambda}}d\mu_{-1}\left(y\right)\label{eq:10}\\
 & =\frac{2\sum_{a=0}^{d-1}\left(-1\right)^{a}\chi\left(a\right)\left(1+\lambda t\right)^{\frac{a}{\lambda}}}{\left(1+\lambda t\right)^{\frac{d}{\lambda}}+1}\left(1+\lambda t\right)^{\frac{x}{\lambda}},\nonumber 
\end{align}
where $d\in\mathbb{N}$ with $d\equiv1\pmod{2}$. 

In view of (\ref{eq:5}), we define the generalized degenerate
Euler polynomials attached to $\chi$ as follows: 
\begin{equation}
\frac{2\sum_{a=0}^{d-1}\left(-1\right)^{a}\chi\left(a\right)\left(1+\lambda t\right)^{\frac{a}{\lambda}}}{\left(1+\lambda t\right)^{\frac{d}{\lambda}}+1}\left(1+\lambda t\right)^{\frac{x}{\lambda}}=\sum_{n=0}^{\infty}\mathcal{E}_{n,\lambda,\chi}\left(x\right)\frac{t^{n}}{n!}.\label{eq:11}
\end{equation}
When $x=0$, $\mathcal{E}_{n,\lambda,\chi}=\mathcal{E}_{n,\lambda,\chi}\left(0\right)$
are called the generalized degenerate Euler numbers attached to $\chi$. 

Let $n$ be an odd natural number. Then, by (\ref{eq:2}), we get
\begin{align}
 & \int_{X}\chi\left(x\right)\left(1+\lambda t\right)^{\frac{nd+x}{\lambda}}d\mu_{-1}\left(x\right)+\int_{X}\chi\left(x\right)\left(1+\lambda t\right)^{\frac{x}{\lambda}}d\mu_{-1}\left(x\right)\label{eq:12}\\
 & =2\sum_{l=0}^{nd-1}\left(-1\right)^{l}\chi\left(l\right)\left(1+\lambda t\right)^{\frac{l}{\lambda}}.\nonumber 
\end{align}

Now, we set 
\begin{equation}
R_{k}\left(n,\lambda\mid x\right)=2\sum_{l=0}^{n}\left(-1\right)^{l}\chi\left(l\right)\left(l\mid\lambda\right)_{k}.\label{eq:13}
\end{equation}

From (\ref{eq:2}) and (\ref{eq:12}), we have 
\begin{align}
 & \int_{X}\left(1+\lambda t\right)^{\frac{x+dn}{\lambda}}\chi\left(x\right)d\mu_{-1}\left(x\right)+\int_{X}\chi\left(x\right)\left(1+\lambda t\right)^{\frac{x}{\lambda}}d\mu_{-1}\left(x\right)\label{eq:14}\\
 & =\frac{2\int_{X}\left(1+\lambda t\right)^{\frac{x}{\lambda}}\chi\left(x\right)d\mu_{-1}\left(x\right)}{\int_{X}\left(1+\lambda t\right)^{\frac{ndx}{\lambda}}d\mu_{-1}\left(x\right)}\nonumber \\
 & =\sum_{k=0}^{\infty}R_{k}\left(nd-1,\lambda\mid\chi\right)\frac{t^{k}}{k!},\nonumber 
\end{align}
where $n,d\in\mathbb{N}$ with $n\equiv1\pmod{2}$, $d\equiv1\pmod{2}$. 

In this paper, we give some identities of symmetry for the generalized
degenerate Euler polynomials attached to $\chi$ derived
from the symmetric properties of certain fermionic $p$-adic integrals
on $\mathbb{Z}_{p}$.

\section{Identities of symmetry for the generalized degenerate Euler polynomials}

Let $w_{1},w_{2}$ be an odd natural numbers. Then we consider the
following integral equation: 
\begin{align}
 & \frac{\int_{X}\int_{X}\left(1+\lambda t\right)^{\frac{w_{1}x_{1}+w_{2}x_{2}}{\lambda}}\chi\left(x_{1}\right)\chi\left(x_{2}\right)d\mu_{-1}\left(x_{1}\right)d\mu_{-1}\left(x_{2}\right)}{\int_{X}\left(1+\lambda t\right)^{\frac{dw_{1}w_{2}x}{\lambda}}d\mu_{-1}\left(x\right)}\label{eq:15}\\
 & =\frac{2\left(\left(1+\lambda t\right)^{\frac{dw_{1}w_{2}}{\lambda}}+1\right)}{\left(\left(1+\lambda t\right)^{\frac{w_{1}d}{\lambda}}+1\right)\left(\left(1+\lambda t\right)^{\frac{w_{2}d}{\lambda}}+1\right)}\nonumber \\
 & \relphantom =\times\sum_{a=0}^{d-1}\chi\left(a\right)\left(1+\lambda t\right)^{\frac{w_{1}a}{\lambda}}\left(-1\right)^{a}\nonumber \\
 & \relphantom =\times\sum_{b=0}^{d-1}\chi\left(b\right)\left(1+\lambda t\right)^{\frac{w_{2}b}{\lambda}}\left(-1\right)^{b}.\nonumber 
\end{align}

From (\ref{eq:10}) and (\ref{eq:11}), we note that 
\begin{equation}
\int_{X}\chi\left(y\right)\left(x+y\mid\lambda\right)_{n}d\mu_{-1}\left(y\right)=\mathcal{E}_{n,\lambda,\chi}\left(x\right),\quad\left(n\ge0\right).\label{eq:16}
\end{equation}

By (\ref{eq:14}), we get 
\begin{equation}
\int_{X}\chi\left(x\right)\left(x+dn\mid\lambda\right)_{k}d\mu_{-1}\left(x\right)+\int_{X}\chi\left(x\right)\left(x\mid\lambda\right)_{k}d\mu_{-1}\left(x\right)=R_{k}\left(nd-1,\lambda\mid x\right),\label{eq:17}
\end{equation}
where $k\ge0$.

Thus, by (\ref{eq:16}) and (\ref{eq:17}), we get 
\begin{equation}
\mathcal{E}_{k,\lambda,\chi}\left(nd\right)+\mathcal{E}_{k,\lambda,\chi}=R_{k}\left(nd-1,\lambda\mid\chi\right),\label{eq:18}
\end{equation}
where $k\ge0$, $n,d\in\mathbb{N}$ with $n\equiv1\pmod{2}$ , $d\equiv1\pmod{2}$. 

Now, we set 
\begin{equation}
I_{\chi}\left(w_{1},w_{2}\mid\lambda\right)=\frac{\int_{X}\int_{X}\chi\left(x_{1}\right)\chi\left(x_{2}\right)\left(1+\lambda t\right)^{\frac{w_{1}x_{1}+w_{2}x_{2}+w_{1}w_{2}x}{\lambda}}d\mu_{-1}\left(x_{1}\right)d\mu_{-1}\left(x_{2}\right)}{\int_{X}\left(1+\lambda t\right)^{\frac{dw_{1}w_{2}x}{\lambda}}d\mu_{-1}\left(x\right)}.\label{eq:19}
\end{equation}

From (\ref{eq:19}), we have 
\begin{align}
 & I_{\chi}\left(w_{1},w_{2}\mid\lambda\right)\label{eq:20}\\
 & =\frac{2\left(\left(1+\lambda t\right)^{\frac{dw_{1}w_{2}}{\lambda}}+1\right)\left(1+\lambda t\right)^{\frac{w_{1}w_{2}x}{\lambda}}}{\left(\left(1+\lambda t\right)^{\frac{w_{1}d}{\lambda}}+1\right)\left(\left(1+\lambda t\right)^{\frac{w_{2}d}{\lambda}}+1\right)}\nonumber \\
 & \relphantom =\times\sum_{a=0}^{d-1}\chi\left(a\right)\left(-1\right)^{a}\left(1+\lambda t\right)^{\frac{w_{1}a}{\lambda}}\nonumber \\
 & \relphantom =\times\sum_{b=0}^{d-1}\chi\left(b\right)\left(-1\right)^{b}\left(1+\lambda t\right)^{\frac{w_{2}b}{\lambda}}.\nonumber 
\end{align}

Thus, by (\ref{eq:20}), we see that $I_{\chi}\left(w_{1},w_{2}\mid\lambda\right)$
is symmetric in $w_{1},w_{2}$. By (\ref{eq:12}), (\ref{eq:14}),
(\ref{eq:16}) and (\ref{eq:19}), we get 
\begin{align}
 & 2I_{\chi}\left(w_{1},w_{2}\mid\lambda\right)\label{eq:21}\\
 & =\sum_{l=0}^{\infty}\left(\sum_{i=0}^{l}\binom{l}{i}\mathcal{E}_{i,\frac{\lambda}{w_{2}},\chi}\left(w_{1}x\right)w_{2}^{i}w_{1}^{l-i}R\left(\left.dw_{2}-1,\frac{\lambda}{w_{1}}\right|\chi\right)\right)\frac{t^{l}}{l!}.\nonumber 
\end{align}

From the symmetric property of $I_{\chi}\left(w_{1},w_{2}\mid\lambda\right)$
in $w_{1}$ and $w_{2}$, we have 
\begin{align}
 & 2I_{\chi}\left(w_{1},w_{2}\mid\lambda\right)\label{eq:22}\\
 & =2I_{\chi}\left(w_{2},w_{1}\mid\chi\right)\nonumber \\
 & =\sum_{l=0}^{\infty}\left(\sum_{i=0}^{l}\binom{l}{i}\mathcal{E}_{i,\frac{\lambda}{w_{1}},\chi}\left(w_{2}x\right)w_{1}^{i}w_{2}^{l-i}R\left(\left.dw_{1}-1,\frac{\lambda}{w_{2}}\right|\chi\right)\right)\frac{t^{l}}{l!}.\nonumber 
\end{align}

Therefore, by (\ref{eq:21}) and (\ref{eq:22}), we obtain the following
theorem.
\begin{thm}
\label{thm:1} For $w_{1},w_{2},d\in\mathbb{N}$ with $w_{1}\equiv w_{2}\equiv d\equiv1\pmod{2}$,
let $\chi$ be a Dirichlet character with conductor $d$. Then,
we have 
\begin{align*}
 & \sum_{i=0}^{l}\binom{l}{i}\mathcal{E}_{i,\frac{\lambda}{w_{1}},\chi}\left(w_{2}x\right)w_{1}^{i}w_{2}^{l-i}R\left(\left.dw_{1}-1,\frac{\lambda}{w_{2}}\right|\chi\right)\\
 & =\sum_{i=0}^{l}\binom{l}{i}\mathcal{E}_{i,\frac{\lambda}{w_{2}},\chi}\left(w_{1}x\right)w_{2}^{i}w_{1}^{l-i}R\left(\left.dw_{2}-1,\frac{\lambda}{w_{1}}\right|\chi\right),
\end{align*}
where $l\ge0$.
\end{thm}
When $x=0$, by Theorem \ref{thm:1}, we get
\begin{align*}
 & \sum_{i=0}^{l}\binom{l}{i}\mathcal{E}_{i,\frac{\lambda}{w_{1}},\chi}w_{1}^{i}w_{2}^{l-i}R\left(\left.dw_{1}-1,\frac{\lambda}{w_{2}}\right|\chi\right)\\
 & =\sum_{i=0}^{l}\binom{l}{i}\mathcal{E}_{i,\frac{\lambda}{w_{2}},\chi}w_{2}^{i}w_{1}^{l-i}R\left(\left.dw_{2}-1,\frac{\lambda}{w_{1}}\right|\chi\right),\quad\left(l\ge0\right).
\end{align*}
By (\ref{eq:19}), we get 
\begin{align}
 & 2I_{\chi}\left(w_{1},w_{2}\mid\lambda\right)\label{eq:23}\\
 & =\sum_{l=0}^{dw_{2}-1}\left(-1\right)^{l}\chi\left(l\right)\int_{X}\left(1+\lambda t\right)^{\frac{w_{2}}{\lambda}\left(w_{2}+w_{1}x+\frac{w_{1}}{w_{2}}l\right)}\chi\left(x_{2}\right)d\mu_{-1}\left(x\right)\nonumber \\
 & =\sum_{n=0}^{\infty}\left(\sum_{l=0}^{dw_{2}-1}\left(-1\right)^{l}\chi\left(l\right)\mathcal{E}_{n,\frac{\lambda}{w_{2}},\chi}\left(w_{1}x+\frac{w_{1}}{w_{2}}l\right)w_{2}^{n}\right)\frac{t^{n}}{n!}.\nonumber 
\end{align}

On the other hand, 
\begin{align}
 & 2I_{\chi}\left(w_{2},w_{1}\mid\lambda\right)=2I_{\chi}\left(w_{1},w_{2}\mid\lambda\right)\label{eq:24}\\
 & =\sum_{n=0}^{\infty}\left(\sum_{l=0}^{dw_{1}-1}\left(-1\right)^{l}\chi\left(l\right)\mathcal{E}_{n,\frac{\lambda}{w_{1}},\chi}\left(w_{2}x+\frac{w_{2}}{w_{1}}l\right)w_{1}^{n}\right)\frac{t^{n}}{n!}.\nonumber 
\end{align}

Therefore, by (\ref{eq:23}) and (\ref{eq:24}), we obtain the following
theorem.
\begin{thm}
\label{thm:2} For $w_{1},w_{2},d\in\mathbb{N}$ with $d\equiv1\pmod{2}$,
$w_{1}\equiv1\pmod{2}$ and $w_{2}\equiv1\pmod{2}$, let $\chi$ be
a Dirichlet character with conductor $d$. Then, we have 
\begin{align*}
 & w_{2}^{n}\sum_{l=0}^{dw_{2}-1}\left(-1\right)^{l}\chi\left(l\right)\mathcal{E}_{n,\frac{\lambda}{w_{2}},\chi}\left(w_{1}x+\frac{w_{1}}{w_{2}}l\right)\\
 & =w_{1}^{n}\sum_{l=0}^{dw_{1}-1}\left(-1\right)^{l}\chi\left(l\right)\mathcal{E}_{n,\frac{\lambda}{w_{1}},\chi}\left(w_{2}x+\frac{w_{2}}{w_{1}}l\right),\quad\left(n\ge0\right).
\end{align*}

\end{thm}
To derive some interesting identities of symmetry for the generalized
degenerate Euler polynomials attached to $\chi$, we used the symmetric properties
for certain fermionic $p$-adic integrals on $\Zp$. When $w_{2}=1$,
from Theorem \ref{thm:2}, we have 
\begin{align*}
 & \sum_{l=0}^{d-1}\left(-1\right)^{l}\chi\left(l\right)\mathcal{E}_{n,\lambda,\chi}\left(w_{1}x+w_{1}l\right)\\
 & =w_{1}^{n}\sum_{l=0}^{dw_{1}-1}\left(-1\right)^{l}\chi\left(l\right)\mathcal{E}_{n,\frac{\lambda}{w_{1}},\chi}\left(x+\frac{1}{w_{1}}l\right).
\end{align*}
In particular, for $x=0$, we get 
\begin{align*}
 & \sum_{l=0}^{d-1}\left(-1\right)^{l}\chi\left(l\right)\mathcal{E}_{n,\lambda,\chi}\left(w_{1}l\right)\\
 & =w_{1}^{n}\sum_{l=0}^{dw_{1}-1}\left(-1\right)^{l}\chi\left(l\right)\mathcal{E}_{n,\frac{\lambda}{w_{1}},\chi}\left(\frac{1}{w_{1}}l\right).
\end{align*}

\bibliographystyle{amsplain}

\providecommand{\bysame}{\leavevmode\hbox to3em{\hrulefill}\thinspace}
\providecommand{\MR}{\relax\ifhmode\unskip\space\fi MR }
\providecommand{\MRhref}[2]{%
  \href{http://www.ams.org/mathscinet-getitem?mr=#1}{#2}
}
\providecommand{\href}[2]{#2}

\end{document}